\documentclass[12pt]{amsart}

\usepackage{graphicx}
\usepackage{stmaryrd}
\usepackage{amssymb}

\usepackage{amsmath,amscd}

\usepackage[T1]{fontenc}
\usepackage{yfonts}

\newcommand{\R}{\mathbb{R}}
\newcommand{\Q}{\mathbb{Q}}

\newcommand{\M}{\mathcal{M}}

\renewcommand{\P}{\mathbb{P}}
\newcommand{\W}{\mathbb{W}}
\renewcommand{\S}{\mathbb{S}}
\renewcommand{\Q}{\mathbb{Q}}
\renewcommand{\M}{\mathbb{M}}
\newcommand{\D}{\mathcal{D}}
\newcommand{\C}{\mathbb{C}}

\newcommand{\SP}{\mathbb{SP}}

\newcommand{\spc}{(\omega+1)^\omega}
\newcommand{\tr}{(\omega+1)^{<\omega}}

\newcommand{\har}{\!\!\upharpoonright\!\!}
\newcommand{\Bor}{\mbox{Bor}}
\newcommand{\dom}{\mbox{dom}}
\newcommand{\diam}{\mbox{diam}}
\newcommand{\dist}{\mbox{dist}}
\newcommand{\cov}{\mbox{cov}}
\renewcommand{\oplus}{(\omega+1)}

\newtheorem{theorem}{Theorem}
\newtheorem{lemma}{Lemma}
\newtheorem{corollary}{Corollary}

\newtheorem{proposition}{Proposition}
\newtheorem*{claim}{Claim}

\theoremstyle{definition}
\newtheorem{definition}{Definition}

\newtheorem{remark}{Remark}
\newtheorem*{example}{Example}

\author{Marcin Sabok} 

\address{Mathematical Institute,
  Wroc\l aw University, pl. Grunwaldzki $2\slash 4$,
  $50$-$384$ Wroc\l aw, Poland }

\email{sabok@math.uni.wroc.pl}

\title{$\sigma$-continuity and related forcings}

\begin{document}

\begin{abstract}
  The Stepr\=ans forcing notion arises as a quotient of
  Borel sets modulo the ideal of $\sigma$-continuity of a
  certain Borel not $\sigma$-continuous function. We give a
  characterization of this forcing in the language of trees
  and using this characterization we establish such
  properties of the forcing as fusion and continuous reading
  of names. Although the latter property is usually implied
  by the fact that the associated ideal is generated by
  closed sets, we show it is not the case with Stepr\=ans
  forcing. We also establish a connection between Stepr\=ans
  forcing and Miller forcing thus giving a new description
  of the latter.  Eventually, we exhibit a variety of
  forcing notions which do not have continuous reading of
  names in any presentation.
\end{abstract}

\maketitle

\section{Introduction}

Many classical forcing notions arise as quotient Boolean
algebras of $\Bor(X)$ modulo an ideal $I$ in a Polish space
$X$. Forcings of this form are called idealized forcings
(cf. \cite{Zpl:FI}) and are usually denoted by $\P_I$ to
indicate the ideal they arise from. In this way Cohen
forcing is associated with the ideal of meager sets, Sacks
forcing with countable sets and Miller forcing with
$K_\sigma$ sets in the Baire space, to recall just a few
examples. The generic extensions given by these forcings are
always extensions by a~single real, which is called the
generic real.

Idealized forcings $\P_I$ often are equivalent to forcings
with certain families of trees ordered by inclusion. For
instance, Sacks forcing is the forcing with perfect trees
and Miller forcing with superperfect trees in
$\omega^{<\omega}$.

In examining forcing effects on the real line it is often
convenient to have a nice representation of names for reals
in the extension. In many examples such a representation is
given by functions from the ground model. Namely each real
in the extension is the value of a certain function from the
ground model at the generic real. Not always, however, can
the function be defined globally. If we assume properness of
the forcing $\P_I$ then we are provided by a representation
in terms of Borel functions:

\begin{theorem}[Zapletal, \cite{Zpl:FI}]
  If the forcing notion $\P_I$ is proper and $\dot x$ is a
  name for a real then for each $B\in \P_I$ there is a
  condition $C\leq B$ and a Borel function
  $f:C\rightarrow\R$ such that $$C\Vdash \dot x=f(\dot g)$$
  where $\dot g$ is the name for generic real.
\end{theorem}

The most desirable situation is when the function can be
chosen to be continuous and in many cases it actually
happens.  This property is called continuous reading of
names. One should be aware, however, that this property
depends (at least formally) on the topology of the space
$X$. How common this property is among idealized forcings,
can be partially accounted for by the following theorem.

\begin{theorem}[Zapletal, \cite{Zpl:FI}]
  If the ideal $I$ is generated by closed sets then the
  associated forcing $\P_I$ is proper and has continuous
  reading of names.
\end{theorem}

There is one important example of a forcing notion $\P_I$
which is proper but fails to have continuous reading of
names in the natural topology of the space. Let us recall
the old problem of Lusin whether there is a Borel function
which is not $\sigma$-continuous. In \cite{CMPS} a
particularly simple example of such a function was given,
namely the Pawlikowski's function $P$. The ideal $I_P$ of
sets on which $P$ is $\sigma$-continuous gives rise to the
forcing notion $\P_{I_P}$, usually called (cf.
\cite{Zpl:FI}) the Stepr\=ans forcing. In \cite{Stpr}
Stepr\=ans introduced this forcing and used it to increase
the cardinal characteristic $\cov(I_P)$ in a generic
extension. The key feature of the forcing $\P_{I_P}$ is that
it adds a real which is not contained in any ground model
set from $I_P$.

Since the ideal $I_P$ can be seen as a porosity ideal (cf.
\cite{Zpl:FI}), properness of the forcing follows from
another general result.

\begin{theorem}[Zapletal, \cite{Zpl:FI}]
  If $I$ is a porosity ideal then the forcing $\P_I$ is
  proper.
\end{theorem}

Stepr\=ans forcing has many nice properties, one of them is
the fact that compact sets are dense in it. This follows
from the following theorem.

\begin{theorem}[Zapletal, \cite{Zpl:DSTDF}]\label{cldense}
  For any Borel not $\sigma$-continuous function
  $f:\omega^\omega\rightarrow\omega^\omega$ and for any
  Borel set $B\not\in I_f$ there exists a compact set
  $C\subseteq B$ such that $C\not\in I_f$.
\end{theorem}

The proof of this theorem introduces a certain Borel game
which detects $\sigma$-continuity of a given Borel function.
The result follows then from determinacy of this game.

Stepr\=ans forcing, however, does not have continuous
reading of names, for the function $P$, treated as a name
for a real, is itself a counterexample.  This single
obstacle may be handled by extending the topology to one
which has the same Borel sets and makes $P$ continuous. A
question is if this would bring about continuous reading of
names in the Stepr\=ans forcing.

This has been investigated by the authors of \cite{Hru.Zpl}
who argued that the ideal associated with Stepr\=ans forcing
is generated by closed sets in the extended topology.  This
should result in continuous reading of names but the
argument from \cite{Hru.Zpl} is incorrect. The problem
whether the ideal is generated by closed sets in this
extended topology was also raised in \cite{Stpr}.

We will show that the ideal $I_P$ is not generated by closed
sets in the extended topology. Nevertheless, we will prove
that the Stepr\=ans forcing has continuous reading of names
in this topology. To this end we will establish a
description of the forcing in terms of trees and deduce
continuous reading of names from properties of these trees.
This will also enable us to define fusion in Stepr\=ans
forcing.

In light of the above another question arises. Are there any
forcing notions of the form $\P_I$ which do not have
continuous reading of names in any presentation (i.e.  in
any Polish topology which gives the same Borel structure)?
This question was already posed in \cite{Hru.Zpl}. Recently
an example has been given by Zapletal in \cite{Zpl:FI},
namely he proved the eventually different real forcing has
this property. We will present a different example.

In fact, we will show that such forcings are quite common
among the idealized forcings. Namely, there is a method of
constructing them out of forcings which, as the Stepr\=ans
forcing, do not have continuous reading of names in one
topology.

\section{Definitions and notation}\label{sec:def}

Throughout this paper an ideal will always mean a
$\sigma$-ideal of subsets of a Polish space.

In a space $X$ a system of sets indexed by a tree
$T\subseteq Y^{<\omega}$ ($Y$ is an arbitrary set). is to be
understood as a map $T\ni\tau\mapsto D_\tau\subseteq X$ such
that if $\tau\subseteq\tau'\in T$ then $D_{\tau'}\subseteq
D_\tau$.  The system is disjoint if $D_\tau\cap
D_{\tau'}=\emptyset$ for $\tau\not=\tau', |\tau|=|\tau'|$.

In a space $X^\omega$, whatever be its topology, for a
finite partial function $\tau:\omega\rightarrow
\mathcal{P}(X)$ we will denote by $[\tau]$ the set $\{t\in
X^\omega: \forall n\in\dom(\tau)\ t(n)\in\tau(n)\}$. For a
tree $T\subseteq X^{<\omega}$ let its limit, denoted $\lim
T$, be the set $\{x\in X^\omega: \forall n\in\omega\quad
x\har n\in T \}$. For a node $\tau\in T$ the end-extension
of $T$ above $\tau$ will stand for the subtree $\{\sigma\in
T: \tau\subseteq\sigma\vee\sigma\subseteq\tau\}$. For a set
$T_0\subseteq T$ the end-extension of $T$ above $T_0$ will
be the union of end-extensions of $T$ above elements of
$T_0$.

We will say that a Borel function $f:X\rightarrow Y$, where
$X,Y$ are Polish spaces, is $\sigma$-continuous if there
exist a countable cover of the space $X=\bigcup_n X_n$ (with
arbitrary sets $X_n$) such that $f\har X_n$ is continuous
for each $n$. It follows from Kuratowski's extension theorem
that we may require that the sets $X_n$ be Borel. If they
can be chosen closed in $X$ then we shall say that $f$ is
closed-$\sigma$-continuous.

If $f:X\rightarrow Y$ is not $\sigma$-continuous then $I_f$
denotes the ideal of sets on which $f$ is
$\sigma$-continuous.

For two Borel functions $f:X\rightarrow Y$ and
$f':X'\rightarrow Y'$ we will say that $f$ can be factorized
by $f'$ if there exist a continuous $1$-$1$ function
$\varphi:X'\rightarrow X$ and an open $1$-$1$
$\psi:Y'\rightarrow Y$ such that the following diagram commutes:
$$
\begin{CD}
  Y' @ >\psi >> Y\\
  @AAf' A @AAfA\\
  X' @> \varphi >> X
\end{CD}
$$

Obviously, if $f'$ is not $\sigma$-continuous and factorizes
$f$ then $f$ is also not $\sigma$-continuous.

In a metric space $(X,d)$ for $A,B\subseteq X$ we will
denote by $\dist(A,B)$ the infimum of $d(a,b)$ for $a\in A$
and $b\in B$. The Hausdorff distance between $A$ and $B$
will be denoted by $h(A,B)$.

The space $\spc$ is endowed with the product topology of
order topologies on $\omega+1$. It is of course homeomorphic
to the Cantor space. We also fix a metric $\rho$ on $\spc$
which gives the above topology. For $x,y\in\spc$ let
$\rho(x,y)=\sum_n \frac{1}{2^n}\rho'(x(n),y(n))$ where
$\rho'$ metrizes $\omega+1$ with its order-topology, i.e.
$\rho'(n,\omega)=\frac{1}{2^n}$ and
$\rho'(n,m)=|\frac{1}{2^n}-\frac{1}{2^m}|$ for $n,m<\omega$.
All metric notions on $\spc$, like diameter, distance, etc.,
will be relative to the metric $\rho$.

The Pawlikowski's function
$P:\spc\rightarrow\omega^\omega$ is defined as follows:
\begin{displaymath}
  P(x)(n) = \left\{ 
    \begin{array}{ll}
      x(n)+1 & \mbox{if}\quad x(n)<\omega,\\
      0 & \mbox{if}\quad x(n)=\omega.
    \end{array} \right.
\end{displaymath}
It has been shown in \cite{CMPS} that $P$ is not
$\sigma$-continuous. Hence $I_P=\{A\subseteq\spc: P\har
A\mbox{ is }\sigma\mbox{-continuous}\}$ is a proper ideal.
It's subideal $I_P^c$ is defined analogously for
closed-$\sigma$-continuity.

Note that the smallest topology on $\spc$ in which $P$ is
continuous is the one with basic clopens of the form
$[\sigma]$ for $\sigma\in(\omega+1)^{<\omega}$. With this
topology $\spc$ is homeomorphic to the Baire space
$\omega^\omega$. We will thus refer to the two topologies:
the original and the extended one as Cantor and Baire
topology, respectively.

Once we have an ideal in $\Bor(X)$ (the family of all Borel
sets in $X$), we consider the associated forcing notion
$\P_I$ which is the poset $(\Bor(X)\setminus I,\subseteq)$.
Of course, it is equivalent to the Boolean algebra
$\Bor(X)\slash I$.  The Stepr\=ans forcing, associated with
$I_P$ will be denoted by $\S$ and the forcing associated
with $I_P^c$ will be denoted by $\S_c$.

We will say that a forcing $\P_I$ has continuous reading of
names in a topology $\mathcal{T}$ of the space $X$ if for
any $B\in\P_I$ and any Borel function
$f:B\rightarrow\omega^\omega$ there exists $\P_I\ni
C\subseteq B$ such that $f\har C$ is continuous in
$\mathcal{T}$.

The general definition of properness of a forcing notion can
be found for instance in \cite{Jech}. Let us, however,
recall a characterization formulated by Zapletal in
\cite{Zpl:FI}: forcing of the form $\P_I$ is proper iff for
every countable elementary substructure $M$ of a large
enough $H_\kappa$ and every condition $B\in M\cap \P_I$ the
set $\{x\in B: x \mbox{ is } \P_I\mbox{-generic over } M\}$
is not in $I$ (this set turns out to be always Borel).

Recall also that a forcing notion $\P$ satisfies
Baumgartner's Axiom A if there is a sequence
$\leq_n,n<\omega$ of partial orders on $\P$ such that
$\leq_0\,=\,\leq$, $\leq_{n+1}\,\subseteq\,\leq_n$ and
\begin{itemize}
\item if $\P\ni p_n, n<\omega$ are such that $p_{n+1}\leq_n
  p_n$ there is a $q\in\P$ such that $q\leq_n p_n$ for all
  $n$,
\item for every $p\in\P$, for every $n$ and for every
  ordinal name $\dot x$ there exist $\P\ni q\leq_n p$ and a
  countable set $B$ such that $q\Vdash \dot x \in B$.
\end{itemize}

\noindent Of course, forcings satisfying Axiom A are proper.

An ideal $I$ in a Polish space $X$ is said to be generated
by closed sets if any Borel $B\in I$ has a $F_\sigma$
superset $C\in I$.

\section{Characterization of the Stepr\=ans forcing}

Note first of all that for every closed set $C$ in the
Cantor topology of $\spc$ there exists a tree
$T\subseteq\tr$ such that $C=\lim T$. Not for every tree,
however, its limit is a closed set in the Cantor topology.
In general, if $T\subseteq(\omega+1)^{<\omega}$ is a tree
then $\lim T$ is a $G_\delta$ set in the Cantor topology
(since its complement is the union of sets $[\tau]$ for
$\tau\not\in T$).

One very popular belief concerning the Stepr\=ans forcing is
that a Borel set $A$ is in $\P_I$ iff it contains limit of a
tree $T\subseteq (\omega+1)^{<\omega}$ such that every
$\tau\in S$ has an extension $\tau'\in S$ which splits into
infinitely many immediate successors including
$\tau'^\smallfrown\omega$. This is, however, not true, as we
will show in the following example.

\begin{example}

  We will now construct a tree $S$ with the property that
  every $\tau\in S$ has an extension $\tau'\in S$ which
  splits into infinitely many immediate successors including
  $\tau'^\smallfrown\omega$ but $P$ is $\sigma$-continuous
  on $\lim S$. We will build the tree inductively on its
  levels.  For any node $\tau\in S$ we will also define a
  set $A_\tau\subseteq\omega$. We begin with $\emptyset$ and
  put $A_\emptyset=\omega$. Suppose that we have the tree
  $S$ built up to level $k$. Now let each node $\tau$ split
  into $\tau^\smallfrown\omega$ as well as $\tau^\smallfrown
  n$ for $n\in A_\tau$. Define sets $A_{\tau^\smallfrown i}$
  for $i\in A_\tau\cup\{\omega\}$ so that they form a
  partition of $A_\tau$ into infinitely many infinite
  subsets. Note at this point that if $s\in\lim S$ and
  $s(n)<\omega$ then $s\har n$ is uniquely determined.
  \begin{claim}
    The function $P$ is $\sigma$-continuous on $\lim S$.
  \end{claim}
  \begin{proof}
    Define $X_n=\{s\in\lim S: \forall m\geq n\quad
    s(m)=\omega\}$ and $X_\infty = \lim S\setminus\bigcup_n
    X_n$. Note that $P$ is continuous on $X_\infty$. Indeed,
    take any convergent sequence $s_n\rightarrow s$ such
    that $s_n,s\in X_\infty$ and notice that if
    $s(m)<\omega$ then $s_n(m)=s(m)$ implies also $s_n\har
    m=s\har m$. Hence, since $[(m,s(m))]$ is a neighborhood
    of $s$, there exists $m'<\omega$ such that $s_n\har m=
    s\har m$ for $n>m'$.  Thus, if $s$ has infinitely many
    values $<\omega$ then the sequence $s_n$ eventually
    stabilizes on each coordinate.  This shows that also
    $P(s_n)\rightarrow P(s)$. Since all sets $X_n$ are
    countable, $P$ is $\sigma$-continuous on $\lim S$.
  \end{proof}
\end{example}

Throughout the rest of this section all topological notions
concerning the space $\spc$ will be relative to the Cantor
topology.

As we have already mentioned, many forcings of the form
$\P_I$ can be equivalently described as tree forcings. It
turns out that Stepr\=ans forcing has similar description in
terms of subtrees of $\tr$.

In fact, the forcing $\SP$ considered by Stepr\=ans in
\cite{Stpr} is actually a forcing with trees.  It adds an
$\S$-generic but it is not clear that $\SP$ is equivalent to
$\S$. But since $\SP$ can be easily seen as a dense subset
of $\Q$ (see below), the equivalence follows from Theorem
\ref{trees}.

\begin{definition}
  A tree $T\subseteq(\omega+1)^{<\omega}$ will be called
  wide if every node $\tau\in T$ has an extension $\tau'\in
  T$ such that the set $\lim T\cap[\tau']$ is nowhere dense
  in $\lim T\cap[\tau]$. A subset of $\spc$ will be called
  wide if it is the limit of a wide tree.
\end{definition}

Obviously, the node $\tau'$ above must be of the form
$\tau{''}^\smallfrown\omega$ for some $\tau''\supseteq\tau$
which splits in $T$ into infinitely many immediate
successors.

Let us denote by $\Q$ the poset of wide trees ordered by
inclusion.

\begin{theorem}\label{trees}
  The Stepr\=ans forcing is equivalent to the forcing $\Q$.
\end{theorem}

The idea to consider wide sets comes from the proof of the
famous theorem of Solecki.

\begin{theorem}[Solecki, \cite{Sol}]\label{pmin}
  For any Baire class 1 function $f:X\rightarrow Y$,
  where $X,Y$ are Polish spaces, either $f$ is
  $\sigma$-continuous or else there exist topological
  embeddings $\varphi$ and $\psi$ such that the following
  diagram commutes:
  $$
  \begin{CD}
    \omega^\omega @ >\psi >> Y\\
    @AAP A @AAfA\\
    (\omega+1)^\omega @> \varphi >> X
  \end{CD}
  $$
\end{theorem}

\begin{proof}[Proof of theorem \ref{trees}] The assertion
  results from Proposition \ref{char} and Proposition
  \ref{factor} given below.
\end{proof}

\begin{proposition}\label{char}
  Assume $B\subseteq \spc$ is Borel such that $P\har B$ is
  not $\sigma$-continuous. Then there exists a wide tree $T$
  such that $\lim T=D\subseteq B$.
\end{proposition}

\begin{proof}

  Let us begin with a claim.

  \begin{claim}\label{ext}
    Let $E\subseteq \spc$ be closed such that $P\har E$ is
    not continuous.  There exists a sequence of disjoint
    relative clopens $C_n, n<\omega$ each of the form
    $[\tau_n]$ such that $\bigcup_n C_n$ is dense in $E$.
  \end{claim}
  \begin{proof}
    Note that any relative open set in $E$ contains a
    relative clopen of the form $[\sigma]$, where
    $\sigma\in(\omega+1)^{<\omega}$. This follows from the
    fact that any basic open set (in $E$) of the form
    $[\tau_1^\smallfrown [n,\omega]^\smallfrown\tau_2]$
    either has a nonempty (in $E$) clopen subset of the form
    $[\tau_1^\smallfrown m^\smallfrown\tau_2]$ for some
    $m\geq n$ or is equal to the relative clopen
    $[\tau_1^\smallfrown \omega^\smallfrown\tau_2]$.  Take
    thus a maximal antichain of clopens (in $E$) of the form
    $[\sigma]$. This antichain can be taken infinite for $E$
    is limit of a tree which is not finitely-branching
    (otherwise $P$ would be continuous on $E$) and thus we
    may extend an antichain given by infinitely many
    immediate successors (by numbers less than $\omega$) of
    a chosen node.
  \end{proof}

  \begin{definition}
    A fusion system in $\spc$ is a tree
    $T\subseteq(\omega+1)^{<\omega}$ together with a family
    of trees $T_\tau, \tau\in T$ such that
    \begin{itemize}
    \item each $\lim T_\tau$ is closed,
    \item $T_\tau$ has stem $\tau$,
    \item for $\tau\subseteq\tau'\in T$\ \
      $T_{\tau'}\subseteq T_\tau$.
    \end{itemize}
  \end{definition}

  Now we pass to the main proof. By Theorem \ref{cldense} we
  may assume that $B=\lim T_\emptyset$ is closed. We may
  also assume that $P$ is not $\sigma$-continuous on any
  basic clopen set. We will construct a fusion system
  $T\subseteq(\omega+1)^{<\omega}$, $T_\tau,\tau\in T$ (and
  denote $D_\tau=\lim T_\tau$) so that the set $D=\lim T$
  will be wide.

  The construction is carried out inductively (beginning
  with $T_\emptyset$) in such a way that having constructed
  $\tau$ and $T_\tau$ we find infinitely many (pairwise
  incomparable) extensions of $\tau$ and appropriate family
  of subtrees of $T_\tau$.

  Suppose we have constructed a node $\tau$ and a tree
  $T_\tau$. By the above Claim we can find an antichain
  $\tau_n, n<\omega$ of extensions of $\tau$ such that
  $\{D_\tau\cap[\tau_n]: n<\omega\}$ is a maximal antichain
  of relative clopens in $D_\tau$. We put $T_{\tau_n}$ to be
  the end-extension above $\tau_n$ in $T$. Let us look now
  at the closed set $E=D_\tau\setminus\bigcup_n D_{\tau_n}$.
  In case $P$ is $\sigma$-continuous on this set we will
  extend $\tau$ by $\tau_n$'s only and call this extension
  regular.  If, however, $P$ is not $\sigma$-continuous on
  $E$ then let us first shrink it to $E'$ by cutting off all
  relative clopens on which $P$ is $\sigma$-continuous. Then
  take any $\tau_\omega\in(\omega+1)^{<\omega}$ which gives
  a nonempty relative clopen in $E'$ (of length > $|\tau|$).
  Now extend $\tau$ additionally by $\tau_\omega$ as well as
  define $T_{\tau_\omega}$ as the tree of
  $E'\cap[\tau_\omega]$.  The extension of this form will be
  called irregular and we will refer to $\tau_\omega$ as the
  irregular node.

  Once the tree has been constructed let us note that each
  node $\tau$ has an irregular extension $\tau'$. Indeed,
  for otherwise $D_\tau$ would be a union of countably many
  closed sets which we have cut off (on each of which $P$
  was $\sigma$-continuous) and the set $[\tau]\cap\lim T$.
  On the latter set, however, $P$ is continuous, hence we
  would get that $P$ is $\sigma$-continuous on $D_\tau$,
  which is not the case. So what is left is to show that the
  set $D_{\tau'}\cap D$ is nowhere dense in $D_\tau\cap D$.
  But its complement contains the intersection of a sequence
  of unions of relative clopens which occur either in
  regular extensions or in irregular as those non-irregular
  ones.  This is, however, an intersection of a sequence of
  dense open sets and it is dense by the Baire category
  theorem (recall that $D$ is a $G_\delta$). So the image of
  the irregular node is nowhere dense, as claimed.
\end{proof}

\begin{remark}
  The fusion method from the above proof will be further
  used to established Axiom A and continuous reading of
  names. We would like to mention, however, that Proposition
  \ref{char} can be also proved without fusion, using the
  method of Cantor-Bendixson analysis instead: Call a tree
  $T\subseteq\tr$ small if for each $\tau\in T$ the set
  $[\tau^\smallfrown\omega]$ is relatively open in $\lim T$.
  It is easy to see that if $T$ is small then $P$ is
  continuous on $\lim T$. Then use a procedure in the
  fashion of the Cantor-Bendixson analysis to cut off from
  $T$ all nodes (and their extensions) such that the
  end-extension above them is small. Then the remaining tree
  will be obviously wide. It is maybe more clear now that
  the wide set which remains is not in the ideal $I_P$.  But
  we do not need this because Proposition \ref{factor} says
  that actually any wide set is $I_P$-positive.
\end{remark}

\begin{remark}
  Yet another way of proving Proposition \ref{char} is to
  apply Theorem \ref{pmin} together with Theorem
  \ref{cldense} and the fact that $P$ is of Baire class 1.
  But the arguments given above are much simpler than the
  proof of Theorem \ref{pmin} in its full strength.
  Nevertheless, one of the ideas from Solecki's proof of
  Theorem \ref{pmin} will be used in the proof of
  Proposition \ref{factor}.
\end{remark}

\begin{proposition}\label{factor}
  Assume $D\subseteq(\omega+1)^\omega$ is a wide set. Then
  there are topological embeddings $\varphi$ and $\psi$ such
  that the following diagram commutes:
   $$
     \begin{CD}
         \omega^\omega @ >\psi >> P[D]\\
         @AAP A @AAP\upharpoonright DA\\
         (\omega+1)^\omega @> \varphi >> D
     \end{CD}
   $$
\end{proposition}

\begin{proof} Let us suppose $D=\lim T$ is wide in $\spc$.
  Let us say that a subtree $S\subseteq T$ is an end-subtree
  if there is a finite set of nodes of $T$ such that $S$ is
  the end-extension of this set. Let $\C$ be the forcing
  with end-subtrees of $T$ ordered by inclusion. $\C$ is of
  course equivalent to the Cohen forcing. Let $M$ be a
  countable elementary submodel of a large enough $H_\kappa$
  such that $P,\C,D\in M$. Let $\D_n, n<\omega$ enumerate
  all dense subsets of $\C$ in $M$. For a dense set
  $\D\subseteq\C$ let $\D^*$ denote the set of all finite
  unions of elements from $\D$.

  We are going to construct only the embedding
  $\varphi:\spc\rightarrow D$ since it already determines
  the function $\psi$. To this end we will define a disjoint
  system of wide $G_\delta$ sets $D_\tau\subseteq D,
  \tau\in\tr$ given as limits of trees $T_\tau\in\C$, such
  that the branches of the system (i.e. $\{T_\tau:
  \tau\subseteq t\}$ for $t\in\spc$) will generate
  $\C$-generic filters over $M$.

  The trees $T_\tau$ will be constructed by induction on
  $|\tau|$ and will satisfy the following conditions:
  \begin{enumerate}
  \item[(i)] $T_\tau\in \D_{|\tau|}^*$,
  \item[(ii)] $\diam(D_\tau)<1\slash |\tau|$,
  \item[(iii)] for each $n$ the map
    $(\omega+1)^n\ni\tau\mapsto D_\tau$ is $h$-continuous.
  \end{enumerate}

  Notice that (i) and (ii) implies that each branch
  generates a generic filter over $M$. Indeed, because any
  extension to an ultrafilter must be generic by (i) and by
  (ii) there is precisely one such extension, since it is
  determined by an appropriate generic real.

  Note at this point that the sets $D_\tau$ suffice to
  construct a factorization. For $t\in\spc$ we define
  $\varphi(t)$ to be the generic real given by the generic
  filter along the branch $t$.  Thanks to (iii) $\varphi$ is
  continuous.  Because of disjointness of the system,
  $\varphi$ is injective, and hence a topological embedding.
  On the other hand, $\psi$ is open because the system is
  disjoint and $P[D_\tau]$ is open in $P[D]$ (since $T_\tau$
  is an end-subtree). To see that $\psi$ is continuous we
  use genericity: a formula of the form $P(\varphi(t))(m)=n$
  is absolute for transitive models and if it holds for $t$
  which is generic over $M$ then it must be forced by some
  condition in the generic filter.

  Before we go on and construct the sets $D_\tau$, let us
  first endow the space $\oplus^n$ with some additional
  structure which will be used in the construction. Let
  $S^n_k\subseteq(\omega+1)^n$ be the set of points of
  Cantor-Bendixson rank $\geq n-k$ (for $k\leq n$). Besides
  the sets $S^n_k$, let us also define a system of
  projections $\pi^n_k:S^n_k\rightarrow S^n_{k-1}$ for
  $1\leq k\leq n$.  If we embed $\oplus^n$ into the cube
  $[0,1]^n$ via $0\mapsto0$, $n\mapsto 1-1\slash n$ for
  $n>0$ then $S^n_k$ can be viewed a set of $n\choose k$ its
  $k$-dimensional faces of the cube. In this setting each
  $S^n_k$ can be projected orthogonally onto $S^n_{k-1}$.
  This projection, however, is ambigous at some diagonal
  points.  Nevertheless, we may pick one of the possible
  values and in this way define a function $\pi^n_k$. In
  other words, if $\tau\in S^n_k\setminus S^n_{k-1}$ then we
  pick one $i\in n$ such that $\tau(i)$ is maximal value
  less than $\omega$ and define $\pi^n_k(\tau)(i)=\omega$
  and $\pi^n_k(\tau)(j)=\tau(j)$ for $j\not= i$. On
  $S^n_{k-1}$ $\pi^n_k$ is the identity. The key feature of
  these functions is that they are continuous, no matter
  which values we have picked.

  \begin{lemma}
    For each $n$ and $1\leq k\leq n$ the projection
    $\pi^n_k: S^n_k\rightarrow S^n_{k-1}$ is continuous.
  \end{lemma}
  \begin{proof}
    Note that any point in $S^n_k$ except
    $(\omega,\ldots,\omega)$ ($k$ times $\omega$) has a
    neighborhood in which projection is unambigous and hence
    continuous. But it is easy to see that at the point
    $(\omega,\ldots,\omega)$ any projection is continuous.
  \end{proof}

  In the construction we will use the following lemma which
  holds in $M$.

  \begin{lemma}\label{extn} Let $S, S'$ be end-subtrees of
    $T$, $D=\lim S, D'=\lim S'$, $\delta>0$ and $k<\omega$.
    \begin{enumerate}
    \item\label{i} There is a sequence $S_i, i\in\omega+1$
      of subtrees of $S$ such that $S_i\in\C$, the sets $D_i=\lim S_i$ are
      disjoint, the map $\omega+1\ni i\mapsto D_i$ is
      $h$-continuous and for each $i$ $\diam(D_i)<\delta$
      and $D_i\in\D_k^*$.
    \item\label{ii} If $S_i$ are as above then there is a
      sequence $S_i',i\in\omega+1$ of subtrees of $S'$ such
      that $S_i'\in\C$, the sets $D_i'=\lim S_i'$ are
      disjoint, the map $\omega+1\ni i\mapsto D_i'$ is
      $h$-continuous, for each $i$ $D_i'\in\D_k^*$,
      $\diam(D_i')<3\delta$ and $h(D_i, D_i')\leq 3h(D,D')$.
    \end{enumerate}
  \end{lemma}
  \begin{proof}
    (\ref{i}) Let us pick any node $\tau^\omega\in S$ such
    that $|\tau^\omega|>k$, $\lim S\cap[\tau^\omega]$ has
    diameter $<\delta\slash 3$ and is nowhere dense in $\lim
    S$ (let $\tau^\omega=\tau^\smallfrown\omega$). Put
    $S_\omega$ equal to the end-subtree of $S$ above
    $\tau^\omega$. Notice that for any $\varepsilon>0$ there
    exists a finite set $\tau_i, i\leq n$ of extensions of
    $\tau^\omega$ such that
    $\diam(D_\omega\cap[\tau_i])<\varepsilon$ for each
    $i\leq n$ and $h(D_\omega,\bigcup_{i\leq n}
    D_\omega\cap[\tau_i])<\varepsilon$. Using the fact that
    $\lim S\cap[\tau^\omega]$ is nowhere dense in $\lim S$
    we can find nodes $\tau_i'$ (being extensions of nodes
    $\tau^\smallfrown n_i$ for some $n_i<\omega$) such that
    $\diam(D\cap[\tau_i'])<\varepsilon$ and
    $\dist(D\cap[\tau_i'],D_\omega\cap[\tau_i])<\varepsilon$
    hold for each $i\leq n$. Now it easily follows that
    $h(D_\omega, \bigcup_{i\leq n}
    D\cap[\tau_i'])<3\varepsilon$ (so in particular the
    first set has diameter $<\delta$ if $\varepsilon$ is
    small enough). We may of course shrink each
    $D\cap[\tau_i']$ so that it is the limit of a tree in
    $\D_k$. This ends the first part of the proof.

    (\ref{ii}) Let $\gamma=h(D,D')$. First we claim that
    there is a finite set of nodes $\tau_i',i\leq n$ in $S'$
    such that
    \begin{itemize}
    \item $\diam(D'\cap[\tau_i'])<\gamma$ for each $i\leq
      n$,
    \item $h(D_\omega,\bigcup_{i\leq n}
      D'\cap[\tau_i'])<2\gamma$,
    \item $\diam(\bigcup_{i\leq n}
      D'\cap[\tau_i'])<3\delta$.
    \end{itemize} 
    Indeed, if $\gamma<\delta$ then we may first find
    finitely many nodes $\tau_i, i\leq n$ in $S_\omega$ and
    then appropriate nodes $\tau_i', i\leq n$ in $S'$ such
    that $h(D_\omega,\bigcup_{i\leq n}
    D_\omega\cap[\tau_i])<2\gamma$, both
    $\diam(D_\omega\cap[\tau_i]),\diam(D'\cap[\tau_i'])<\delta$
    and also
    $\dist(D_\omega\cap[\tau_i],D'\cap[\tau_i'])<\delta$ for
    $i\leq n$. Then by the triangle inequality
    $\diam(\bigcup_{i\leq n}D'\\ \cap[\tau_i'])<3\delta$.
    If $\gamma\geq\delta$ then we will do by picking in a
    similar manner just one node $\tau$ and $\tau'$ in $S,
    S'$ respectively.
    
    Now above each $\tau_i'$ choose a node
    $\tau^\omega_i{'}$ such that $D'\cap[\tau^\omega_i{'}]$
    is nowhere dense in $D'$.  Then $\bigcup_{i\leq n}
    D'\cap[\tau^\omega_i{'}]$ is also nowhere dense and has
    diameter $<3\delta$. We may now find conditions in
    $\D_k$ which are stronger than the end-extensions in
    $S'$ above the nodes $\tau^\omega_i{'}$. And put
    $S_\omega'$ to be the union of these, $D_\omega'=\lim
    S_\omega'$. It is clear that
    $h(D_\omega,D_\omega')<3\gamma$. Now as in (\ref{i}) we
    can find a sequence of disjoint subtrees $S_n'\in\D_k^*$
    such that if we put $D_n'=\lim S_n'$ then
    $h(D_n',D_\omega)<1\slash n$. Now it follows from the
    triangle inequality that $h(D_n,D_n')<3\gamma$ as well
    as $\diam(D_n')<3\delta$ holds for $n$ big enough.  But
    we may change those finitely many $S_n$'s (in the same
    way we have found $S_\omega'$, possibly shrinking the
    existing sets) to ensure that this holds for all $n$.
  \end{proof}

  Now we proceed as follows. Let $D_\emptyset=D$. Having
  defined $D_\tau$ for all $\tau\in(\omega+1)^n$ we define
  it for $\tau\in(\omega+1)^{n+1}$. This in turn is done by
  another induction on the sets $S_k\times(\omega+1)$ for
  $0\leq k\leq n$.  That is, we first define $D_\tau$ for
  $\tau\in S^n_0\times(\omega+1)$ and then show how to
  extend the definition from $S^n_k\times(\omega+1)$ to
  $S^n_{k+1}\times(\omega+1)$.  During this construction we
  take care that for each $k$ and $\tau\in S^n_k$
  \begin{equation}\label{ineq}
    \diam(D_\tau)<1\slash
    (3^{n-k}(n+1)) 
  \end{equation}
  and the map $S_k\times(\omega+1)\ni\tau\mapsto D_\tau$ is
  $h$-continuous.

  To start with we use Lemma \ref{extn}(\ref{i}). Suppose we
  have $D_\tau$ defined for $\tau\in S_k\times(\omega+1)$.
  Let us abbreviate $\pi^n_k:S^n_{k+1}\rightarrow S^n_k$ by
  $\pi$ for a moment. For each $\tau\in S^n_{k+1}$ we use
  Lemma \ref{extn}(\ref{ii}) for $D_\tau$ and
  $D_{\pi(\tau)}$ to find sets $D_{\tau^\smallfrown i}$ for
  $i\in\omega+1$ such that
  \begin{equation}\label{three}
    h(D_{\tau^\smallfrown i},D_{\pi(\tau)^\smallfrown
      i})\leq 3 h(D_\tau,D_{\pi(\tau)}),
  \end{equation}
  $\diam(D_{\tau^\smallfrown i})<3\,\diam(D_\tau)$ and
  $D_{\tau^\smallfrown i}\in \D_{|\tau|}^*$. Now
  (\ref{ineq}) follows from the inductive assumption. To see
  $h$-continuity notice that if $(\tau_n,i_n)\rightarrow
  (\tau,i)$ is a convergent sequence in
  $S^n_{k+1}\times(\omega+1)$ then either $\tau_n$ is
  eventually constant, in which case the continuity is easy,
  or $\tau\in S^n_k$ and then $\pi(\tau_n)\rightarrow\tau$
  thanks to the continuity of $\pi$. But then the assertion
  follows from the induction assumption, (\ref{three}) and
  the triangle inequality:
  \begin{eqnarray*}
  h(D_{\tau_n^\smallfrown i_n}, D_{\tau^\smallfrown
    i})\leq\!\!\!&h(D_{\tau_n^\smallfrown
    i_n},D_{\pi(\tau_n)^\smallfrown
    i_n})+h(D_{\pi(\tau_n)^\smallfrown
    i_n},D_{\tau^\smallfrown i_n})\\ &+h(D_{\tau^\smallfrown
    i_n},D_{\tau^\smallfrown i}).  
  \end{eqnarray*}
  In this way we have constructed the sets $D_\tau$ and
  finished the proof.
\end{proof}

Propositions \ref{trees} and \ref{factor} have the following
corollaries.

\begin{corollary}
  If $B\subseteq\spc$ is a Borel set such that $P\har B$ is
  not $\sigma$-continuous then there exists a closed wide
  set $D\subseteq B$.
\end{corollary}

\begin{proof}
  This follows from Proposition \ref{factor} since the image
  of $\varphi$ is a wide closed set.
\end{proof}

The second corollary is a particular case of Theorem
\ref{pmin} when $f$ is the restriction of $P$ to a
$I_P$-large set.

\begin{corollary}
  If $D\subseteq\spc$ is Borel then either $P\har D$ is
  $\sigma$-continuous or else there are topological
  embeddings $\varphi$ and $\psi$ such that the following
  diagram commutes:
   $$
     \begin{CD}
         \omega^\omega @ >\psi >> P[D]\\
         @AAP A @AAP\upharpoonright DA\\
         (\omega+1)^\omega @> \varphi >> D
     \end{CD}
   $$
\end{corollary}

\section{Continuous reading of names}

Let us recall now that Stepr\=ans forcing does not have
continuous reading of names in the Cantor topology of
$\spc$. We will now show that it has continous reading of
names in the Baire topology.

\begin{theorem}\label{crn}
  The forcing notion $\S$ has continuous reading of names in
  the Baire topology on $\spc$.
\end{theorem}

Recall that all metric notions on $\spc$ (like diameter,
distance, etc.) are relative to the metric $\rho$ (see
Section \ref{sec:def}) on $\spc$.

\begin{proof}
  Let $B$ be any Borel $I_P$-positive set in $\spc$ and
  $\dot x$ be a $\S$-name for a real. By Proposition
  \ref{char} we may assume $B$ is a limit of a wide tree.
  Continuous reading of names will result from the following
  claim.
  \begin{claim}\label{dist}
    Let $T$ be a tree and $\sigma\in T$ be such that
    $[\sigma^\smallfrown\omega]\cap\lim T$ is nowhere dense
    in $[\sigma]\cap\lim T$. Then for each $\tau\in T$ such
    that $\sigma^\smallfrown\omega\subseteq\tau$, any
    $\varepsilon>0$ and $n<\omega$ there is $m>n$ and
    $\sigma^\smallfrown m\subseteq\tau'\in T$ such that
    \begin{itemize}
    \item $[\tau']\cap\lim T$ is a relative clopen,
    \item $\diam([\tau']\cap\lim T)<\varepsilon\slash 2$,
    \item $\dist([\tau]\cap\lim T,[\tau']\cap\lim
      T)<\varepsilon\slash 2$.
    \end{itemize}
  \end{claim}
  \begin{proof}
    Consider the family of relative clopen sets of the form
    $[\tau']\cap\lim T$ with $\tau'$ extending some
    $\sigma^\smallfrown m, m>n$ and having diameters
    $<\varepsilon\slash 2$. Put also $\delta=\inf_{i\leq
      n}\dist([\tau], [\sigma^\smallfrown i])$. If the
    assertion of this lemma were false, then the open ball
    around $[\tau]$ with radius
    $\min\{\delta,\varepsilon\slash 2\}$ would exhibit that
    $\sigma^\smallfrown\omega$ has nonempty interior.
  \end{proof}
  Now let us finish the proof of the theorem. We will say
  that a set of nodes of a tree $S$ is a spanning set if $S$
  is the smallest tree containing those nodes. We will find
  a wide subtree $S\subseteq T$ and a spanning set
  $\{\tau_\sigma: \sigma\in\omega^\omega\}$ (with
  $\sigma\mapsto\tau_\sigma$ being order isomorphism) of
  nodes of $S$ together with a set of natural numbers
  $k_\tau$ such that for $\tau$ in the spanning set $$(\lim
  S\cap[\tau])\Vdash_\S \dot x(|\tau|)=k_\tau.$$ This will
  clearly show that on $\lim S$ the name $\dot x$ is read by
  a function continuous in the Baire topology relativized to
  $\lim S$.

  The construction is conducted inductively on $|\tau|$ in a
  fusion manner, that is we define additionally subtrees
  $T_\tau$ with stems $\tau$, respectively. Suppose we have
  found everything up to the level $n$. We will show how to
  extend a single node. First find an extension
  $\tau\subseteq\tau'\in T_\tau$ such that
  $[\tau'^\smallfrown\omega]\cap T_\tau$ is nowhere dense in
  $[\tau]\cap T_\tau$. Without loss of generality let us
  assume that $\{n<\omega: \tau'^\smallfrown n\in T_\tau\}$
  is the whole $\omega$. Now find a forcing extension of
  $[\tau'^\smallfrown\omega]\cap T_\tau$ to a limit of a
  wide tree $T_{\tau^\smallfrown 0}$ such that for some
  $k_{\tau^\smallfrown 0}$
  $$\lim T_{\tau^\smallfrown 0}\Vdash_\S\dot x(|\tau|+1)=
  k_{\tau^\smallfrown 0}.$$ Next using the above lemma and
  some bookkeeping find extensions $\tau'^\smallfrown
  n\subseteq \tau_n', n\in\omega$ so that for any $\sigma\in
  T_{\tau^\smallfrown 0}$ and $n<\omega$ there is
  $m\in\omega$ such that $\dist(\lim
  T_\tau\cap[\tau'_m],\lim T_{\tau^\smallfrown
    0}\cap[\sigma])<1\slash n$ and $\diam(\lim
  T_\tau\cap[\tau'_m])<1\slash n$. Then extend the forcing
  conditions $[\tau'_m]\cap \lim T_\tau$ to limits of wide
  trees $T_{\tau^\smallfrown m+1}$ such that for some
  natural numbers $k_{\tau^\smallfrown m+1}$
  $$\lim T_{\tau^\smallfrown m+1}\Vdash_\S\dot x(|\tau|+1)=
  k_{\tau^\smallfrown m+1}.$$ Notice that if $$\diam(\lim
  T_\tau\cap[\tau'_m]),\ \dist(\lim T_\tau\cap[\tau'_m],\lim
  T_{\tau^\smallfrown 0}\cap[\sigma])<1\slash 2n$$ then
  $$\dist(\lim T_{\tau^\smallfrown m+1},\lim
  T_{\tau^\smallfrown 0}\cap[\sigma])<1\slash n,$$ so the
  interior of $T_{\tau^\smallfrown 0}$ remains empty.
  Moreover, it will remain empty even when we pass to the
  fusion tree $S$.  Thus after the fusion we get an
  $I_P$-positive set and numbers $k_\tau$ which define a
  continuous function in the Baire topology.
\end{proof}

\section{The fusion}

It has been established by Zapletal both in \cite{Zpl:FI}
and in \cite{Zpl:DSTDF} that Stepr\=ans forcing is proper.
The fusion method used in Theorem \ref{trees} and Theorem
\ref{crn} suggests, however, that Axiom A can be deduced
quite easily once we have the notion of a wide tree. Axiom A
was also established by Stepr\=ans for the forcing $\SP$
considered in \cite{Stpr}. We present a different proof that
seems more natural for the forcing of wide trees.

\begin{theorem}
  The Stepr\=ans forcing notion satisfies Axiom A.
\end{theorem}
\begin{proof}
  Let $\W'$ be a forcing with trees $T$ satisfying the
  following conditions:
  \begin{enumerate}
  \item each $\tau\in T$ either has only one immediate
    successor or is such that $\tau^\smallfrown\omega\in T$
    and $[\tau^\smallfrown\lim T]\cap\lim T$ is nowhere
    dense in $\lim T$,
  \item whenever $\tau\in T$ is such that
    $\tau^\smallfrown\omega\in T$ and $[\tau^\smallfrown\lim
    T]\cap\lim T$ is nowhere dense in $\lim T$, we have the
    following. For each $n<\omega$ such that
    $\tau^\smallfrown n\in T$ denote the stem of the tree
    $T$ above $\tau^\smallfrown n$ by $\tau_n$. For each
    clopen $C$ intersecting $\lim T\cap
    [\tau^\smallfrown\omega]$ and any $\varepsilon>0$ there
    is $n<\omega$ such that $\diam [\tau_n]<\varepsilon$ and
    $\dist([\tau_n],C)<\varepsilon$.
  \end{enumerate}
  It is easy to see that $\W'$ is dense in $\W$. So it is
  enough to show that $\W'$ satisfies Axiom A.

  Let $\prec$ be a linear order on $(\omega+1)^{<\omega}$
  such that each $\tau$ occurs later than its initial
  segments. For a tree $T\in\W'$ let $w(T)$ be the set of
  those nodes of $T$ which have more than one immediate
  successor. For $\tau\in w(T)$ the set of its
  \textit{immediate successors in $w(T)$} stands for the set
  $\{\tau'\in w(T):\neg\exists\tau''\in w(T)\quad
  \tau\subsetneq\tau''\subsetneq\tau'\}$. Now let us denote
  by $w_n(T)$ the set of $n$ first (with respect to $\prec$)
  elements of $w(T)$ together with theirs immediate
  successors in $w(T)$.
  
  For $T,S\in\W'$ let $T\leq_n S$ if $T\leq S$ and
  $w_n(S)\subseteq T$. It is now easy to see that with these
  orderings $\W'$ satisfies Axiom A.
\end{proof}

\section{Generating by closed sets}

Both properness and continuous reading of names could be
deduced more easily if only we knew that $I_P$ were
generated by closed sets in the Baire topology. A typical
mistake that may lead to such a conclusion is the conviction
that if $A\subseteq\spc$ is such that $P\har A$ is
continuous in the Baire topology then so it is on the
closure of $A$ in the Baire topology. This is not true, as
observed by Pawlikowski (in a private conversation). 

The next proposition says that continuous reading of names
can hold even when the ideal is not generated by closed
sets.

\begin{proposition}\label{clgen}
  The ideal $I_P$ is not generated by closed sets in the
  Baire topology.
\end{proposition}
Throughout this proof let $\overline X$ (for $X\subseteq
(\omega+1)^\omega$) denote the closure of $X$ in the Baire
topology.
\begin{proof}
  Let us first consider the following set $A=\{\alpha_n,\beta_n:n<\omega\}\subseteq\spc$, where
$$\alpha_n(0)=n,\quad\alpha_n(k)=\omega\ \mbox{for}\ k>0$$
$$\beta_n(n)=0,\quad\beta_n(k)=\omega\ \mbox{for}\
k\not=n.$$ Note that $P\har A$ is continuous. On the other
hand, $P\har\overline{A}$ is not continuous since
$\alpha=(\omega,\omega,\ldots)\in\overline{A}$ and
$\alpha_n\rightarrow\alpha$, whereas
$P(\alpha_n)\not\rightarrow P(\alpha)$.

Using a bijection from $\omega$ to $\omega\times\omega$ we
may identify $(\omega+1)^\omega$ with
$(\omega+1)^{\omega\times\omega}\simeq((\omega+1)^\omega)^\omega$.
Under this identification $P$ becomes $\prod_{n<\omega}P$,
which we will denote by $P^\omega$. First note that
$P^\omega$ is continuous on $A^\omega$ as a product of
continuous functions, so $A^\omega\in I_{P^\omega}$. We will
prove, however, that $A^\omega$ cannot be covered by
countably many sets $F_n, n<\omega$ closed in the Baire
topology, with each $F_n\in I_{P^\omega}$.

Suppose that $A^\omega\subseteq \bigcup_n F_n$ and $F_n$ are
closed in the Baire topology. As $A$ is a discrete set in
$(\omega+1)^\omega$ (in both topologies), the relative
topology (with respect to any of these two) on $A^\omega$ is
that of the Baire space.  $F_n\cap A^\omega$ are relatively
closed, so by the Baire category theorem, one of them has
nonempty interior. This means that there is $n<\omega$,
$k<\omega$ and $\alpha\in A^k$ such that $\alpha^\smallfrown
A^{\omega\setminus k}\subseteq F_n$.  Without loss of
generality $k=0$ and $\overline{A^\omega}\subseteq F_n$. But
$\overline{A^\omega}=(\overline A)^\omega$ and $\overline A$
contains a convergent sequence $\alpha_n\rightarrow \alpha$
such that $P(\alpha_n)\not\rightarrow P(\alpha)$. So if
$A'=\{\alpha,\alpha_n: n<\omega\}$ then $P[A']$ is a
discrete set and $(A')^\omega\subseteq F_n$. Notice,
however, that $P^\omega\har (A')^\omega=(P\har A')^\omega$
is not $\sigma$-continuous, since it can obviously be
factorized by $P$. Hence $F_n\not\in I_{P^\omega}$, which
ends the proof.
\end{proof}

\section{Connections with the Miller forcing}

A natural question that arises after realizing that
Stepr\=ans forcing is described in terms of wide trees is
whether this forcing is equivalent to the Miller forcing. A
negative answer follows for instance from Proposition
\ref{addreal} below. It turns out, however, that Miller
forcing is very close to $\sigma$-continuity, namely it is
isomorphic to the forcing associated to the ideal of
closed-$\sigma$-continuity of $P$.

\begin{proposition}
  The forcing notion $\S_c$ is equivalent to the Miller
  forcing notion.
\end{proposition}
\begin{proof}
  The Miller forcing is equivalent to the forcing
  $\Bor(\omega^\omega)\slash K_\sigma$ and $\S_c$ is
  equivalent to $\Bor((\omega+1)^\omega)\slash I^c_P$. The
  isomorphism is given by the function $P$ itself (as it
  gives rise to a Borel isomorphism of the spaces). The only
  thing to realize is that for $A\subseteq(\omega+1)^\omega$
  $P[A]$ is compact if and only if $A$ is a closed set on
  which $P$ is continuous.  But this is the case since a
  continuous image of a compact set is compact, $P^{-1}$ is
  continuous and a continuous bijection defined on a compact
  set is a homeomorphism.
\end{proof}

The following definition is a natural generalization of
well known notions like Cohen real, Miller real, etc.

\begin{definition}
  Let $M$ be a transitive countable model. We say that
  $s\in\spc$ is a Stepr\=ans real over $M$ if $s\not\in A$
  for any $A\subseteq\spc$ such that $A\in I_P$ and $A$ is
  coded in $M$.
\end{definition}

It is obvious that the generic real for Stepr\=ans forcing
is a Stepr\=ans real over the ground model. In order to
distinguish Stepr\=ans forcing from Miller forcing let us
prove that there are no Stepr\=ans reals in extensions by a
single Miller real.

\begin{proposition}\label{addreal}
  Miller forcing does not add Stepr\=ans real.
\end{proposition}

\begin{proof}
  Let us denote the Miller forcing notion by $\M$. Suppose,
  towards a contradiction, that $\dot s$ is a $\M$-name for
  a Stepr\=ans real. Since $\spc\simeq
  2^\omega\subseteq\omega^\omega$, $\dot s$ is a name for an
  element of $\omega^\omega$. Since Miller forcing has
  continuous reading of names we have a forcing condition
  $B\subseteq\omega^\omega$ and a continuous function
  $f:B\rightarrow\omega^\omega$ such that $$D\Vdash\dot s=
  f(\dot m),$$ where $\dot m$ is the name for the
  $\M$-generic real. By another well known property of Miler
  forcing there exists a stronger condition $D\subseteq B$
  such that either $f\har D$ is constant or $f\har D$ is a
  topological embedding. We can exclude the first
  possibility. Let us denote $E=f[D]$ and note that since
  $P\har E$ is Borel, there is a dense $G_\delta$ set
  $G\subseteq E$ such that $P\har G$ is continuous. But then
  $f^{-1}[G]$ is comeager in $D$ and hence is a condition in
  Miller forcing. But $$f^{-1}[G]\Vdash \dot s\in G$$ and
  $G\in I_P$ which gives a contradition.
\end{proof}

\section{A forcing without continuous reading of names in
  any presentation}

The Stepr\=ans forcing notion does not have continuous
reading of names in one presentation and has it in another.
Let us now show how to use Stepr\=ans forcing to produce a
forcing $\P_I$ which is proper and does not have continuous
reading of names in any presentation.

\begin{theorem}
  There exist an ideal $I\subseteq\Bor(\omega^\omega)$ such
  that the forcing $\P_I$ is proper but it does not have
  continuous reading of names in any presentation.
\end{theorem}
\begin{proof}
  First notice that any presentation of a Polish space $X$
  is given by a Borel isomorphism with another Polish space
  $Y$ and the latter can be assumed to be a $G_\delta$
  subset of $[0,1]^\omega$. Instead of $\omega^\omega$ let
  us consider $X=(\omega^\omega)^2$ with its product
  topology.  Note that each $G_\delta$ set in
  $[0,1]^\omega$ as well as a Borel isomorpshism from
  $\omega^\omega$ to $X$ can be coded by a real. Let
  $x\in\omega^\omega$ code a pair $(G_x,f_x)$ that defines a
  presentation of $X$ as above, i.e. $f_x:G_x\rightarrow X$.
  For $x\in\omega^\omega$ $f_x^{-1}[X_x]$ ($X_x$ denotes the
  vertical section of $X$ at $x$) is an uncountable Borel
  set in $G_x$ and contains a copy $C_x$ of $\spc$.  Let
  $I_x$ be the transported ideal $I_P$ from $C_x$ to $X_x$.
  We define an ideal $I$ on $\Bor(X)$ as follows:
  $$I=\{A\in\Bor(X): \forall x\in\omega^\omega\, A_x\in I_x\}.
  $$ $\P_I$ does not have continuous reading of names in any
  presentation for if $(G_x,f_x)$ defines a presentation
  then $(P\circ f_x^{-1})\har (f_x [C_x])$ is a
  counterexample to continuous reading of names in its
  topology. Let us show that $\P_I$ is proper. Take
  $H_\kappa$ big enough so that the function
  $\omega^\omega\ni x\mapsto (G_x,f_x,C_x,I_x)$ is in
  $H_\kappa$. Take any model $M\prec H_\kappa$ which
  contains this function and let $B\in M$ be $I$-positive.
  By elementarity there is $x\in\omega^\omega\cap M$ such
  that $B_x\not\in I_x$. So $B_x\cap f_x^{-1}[C_x]$ is
  $I_x$-positive. Now since the forcing $\P_I$ below
  $f_x^{-1}[C_x]$ is equivalent to the Stepr\=ans forcing,
  which is proper, it follows that there is $B'\subseteq
  B_x\cap f_x^{-1}[C_x]$ $I_x$-positive and $M$-master.
  
\end{proof}

\section{Acknowledgements}

The author would like to thank Janusz Pawlikowski for
introducing into the subject and for continuous
encouragement.


\begin{thebibliography}{}

\bibitem{CMPS} Cicho\'n J., Morayne M., Pawlikowski J. and
  Solecki S., \textit{Decomposing Baire functions}, Journal
  of Symbolic Logic, Vol. 56, 1991
\bibitem{Hru.Zpl} Hru\v s\'ak M. and Zapletal J.,
  \textit{Forcing with quotients},
  \texttt{arXiv:math/0407182v1}
\bibitem{Jech} Jech T., \textit{Set theory. The third
    Millenium edition, revised and expanded}, Springer, 2006
\bibitem{Sol} Solecki S., \textit{Decomposing Borel sets and
    functions and the structure of Baire class~1 functions},
  Journal of the American Mathematical Society, Vol. 11, No.
  3, 1998
\bibitem{Stpr} Stepr\=ans J., \textit{A very discontinuous
    Borel function}, Journal of Symbolic Logic, Vol. 58, No.
  4., 1993
\bibitem{Zpl:DSTDF} Zapletal J., \textit{Descriptive Set
    Theory and Definable Forcing}, Memoirs of the American
  Mathematical Society, 2004
\bibitem{Zpl:FI} Zapletal J., \textit{Forcing Idealized},
  Cambridge Tracts in Mathematics 174, 2008
\end{thebibliography}
\end{document}